\documentclass{amsart}
\usepackage{amssymb}
\usepackage{amsthm}
\usepackage{amscd}
\usepackage{graphics}

\newtheorem*{reflemma}{Lemma}
\newtheorem*{refprop}{Proposition}
\newtheorem{introtheorem}{Theorem}
\newtheorem{theorem}{Theorem}

\newtheorem*{ML}{Main Lemma}
\newtheorem{lemma}{Lemma}[section]
\newtheorem{cor}{Corollary}
\newtheorem{introcor}{Corollary}
\newtheorem{prop}{Proposition}

\numberwithin{equation}{section}
\title{Symplectic embeddings of polydisks}
\author{Larry Guth}
\address{Department of Mathematics, Stanford, Stanford CA, 94305 USA}
\email{lguth@math.stanford.edu}

\begin{document}
\begin{abstract} If $P$ is a polydisk with radii $R_1 \le ... \le
R_n$ and $P'$ is a polydisk with radii $R_1' \le ... \le R_n'$,
then we prove that $P$ symplectically embeds in $P'$ provided
that $C(n) R_1 \le R_1'$ and $C(n) R_1 ... R_n \le R_1' ...
R_n'$.  Up to a constant factor, these conditions are optimal.
\end{abstract}

\maketitle

\section{Introduction}

In this paper, we study when it is possible to symplectically
embed one polydisk into another.  The volume gives one
obstruction to finding an embedding, and Gromov's non-squeezing
theorem \cite{G1} gives a second obstruction. We prove that, up
to a constant factor, these are the only obstructions.

We will work in $\mathbb{R}^{2n}$ with coordinates $x_1, ... x_n$
and $y_1, ..., y_n$, and with the standard symplectic form
$\omega = \sum_{i=1}^n dx_i \wedge dy_i$.  We let $P$ denote the
polydisk $B^2(R_1) \times ... \times B^2(R_n)$, defined by the
inequalities $x_1^2 + y_1^2 < R_1^2$, ..., $x_n^2 + y_n^2 <
R_n^2$.  We assume that the radii are ordered so that $R_1 \le
... \le R_n$.  Similarly, we let $P'$ denote the polydisk
$B^2(R_1') \times ... \times B^2(R_n')$, with $R_1' \le ... \le
R_n'$.  If $P$ symplectically embeds in $P'$, then the
conservation of volume implies $R_1 ... R_n \le R_1' ... R_n'$
and the non-squeezing theorem implies $R_1 \le R_1'$.

\begin{introtheorem} There is a dimensional constant $C(n)$ so
that the following holds.  If $C(n) R_1 \le R_1'$ and $C(n) R_1
... R_n \le R_1' ... R_n'$ then $P$ symplectically embeds into
$P'$.
\end{introtheorem}

There is a linear embedding from $P$ into $P'$ roughly if and
only if $R_k \lesssim R_k'$ for every $k$.  The first interesting
non-linear embeddings were constructed by Traynor in \cite{T}. 
Using Traynor's methods, one can embed $P$ into $P'$ roughly if
and only if $R_1 ... R_k \lesssim R_1' ... R_k'$ for every $k$
between $1$ and $n$.  Until this paper, it was not known whether
the conditions $R_1 ... R_k \lesssim R_1' ... R_k'$ were
necessary for $1 < k < n$.

The problem of determining exactly when one polydisk
symplectically embeds in another looks very difficult.  The known
results are discussed in detail in \cite{H2}, and many other
embedding problems are discussed in \cite{S}.

\vskip5pt

{\bf Intermediate capacities}

\vskip5pt

Symplectic embeddings are closely related to symplectic
capacities.  Modulo details, a (generalized) symplectic capacity
is a function $c$ that assigns a number to each open set in
$\mathbb{R}^{2n}$ in such a way that if $U$ symplectically embeds
in $V$ then $c(U) \le c(V)$.  The function $c$ should also scale
in a reasonable way if we scale the set $U$, and it obeys some
other mild conditions. (See \cite{H2} for all details.)  The
volume is a trivial example of a generalized capacity.  The first
non-trivial capacity comes from Gromov's proof of the
non-squeezing theorem.  In \cite{EH}, Ekeland and Hofer
constructed an infinite sequence of symplectic capacities.  Since
then, several authors have constructed new capacities.  The
review paper \cite{H2} gives a survey of the field.  Finding new
capacities is an important topic in symplectic geometry.  If
Theorem 1 had been false, then there would have been other
generalized capacities, significantly different from the known ones.

Many capacities involve the 2-dimensional area of some object,
such as a pseudo-holomorphic curve.  On the other hand, the
2n-dimensional volume of a region is a generalized capacity.  In
\cite{H}, Hofer asked whether there are intermediate capacities
that involve 2k-dimensional volumes for $1 < k < n$.  Hofer
defined an intermediate capacity of dimension $k$ to be a
generalized capacity $c$ with the following properties:

\vskip3pt

$1.$ \hskip5pt $c[B^{2k}(1) \times \mathbb{R}^{2n-2k}] < \infty$.

$2.$ \hskip5pt $c[B^{2k-2}(1) \times \mathbb{R}^{2n-2k+2}] = \infty$.

\vskip3pt

For $k=1$, these conditions are satisfied by many interesting
examples as described in \cite{H2}.  For $k=n$, these conditions
are satisfied by the volume.  But for the intermediate range,
there are no examples known of such a capacity.  If the capacity
$c$ is also reasonably continuous, then it would satisfy the
following slightly stronger properties.

\vskip3pt

$1'.$ \hskip5pt $\lim_{R \rightarrow \infty} c[B^{2k}(1) \times B^{2n-2k}(R)]
< \infty.$

$2'.$ \hskip5pt $\lim_{R \rightarrow \infty} c[B^{2k-2}(1) \times
B^{2n-2k+2}(R)] = \infty$.

\vskip3pt

If $k$ lies in the intermediate range $1 < k < n$, Theorem 1
implies that there are no reasonably continuous capacities of
dimension $k$.

\vskip5pt

{\bf Expanding embeddings and symplectic embeddings}

\vskip5pt

Let $X$ denote an open set in $\mathbb{R}^n$.  The cotangent
bundle $T^* X$ is a symplectic manifold which we can think of as
the set of all pairs $(x,y)$ with $x \in X$ and $y \in
\mathbb{R}^n$, equipped with the standard symplectic structure
$\omega$.  We write $U^* X$ to denote the unit ball cotangent
bundle, consisting of all pairs $(x,y)$ with $x \in X$ and $|y| <
1$.

Any smooth embedding $I$ from $X$ to $X'$ induces a symplectic
embedding from $T^* X$ into $T^* X'$.  We say that $I$ is an
expanding embedding if $I$ increases the length of every curve in
$X$, or equivalently if $|dI(v)| \ge |v|$ for every tangent
vector $v$ in $TX$.  If $I$ is an expanding embedding, then it
induces a symplectic embedding from $U^* X$ into $U^* X'$.

If there is no expanding embedding from $X$ to $X'$, then it's
interesting to know whether there is a symplectic embedding from
$U^* X$ into $U^* X'$.  If there is no symplectic embedding, then
we can say that the symplectic geometry has remembered that $X$
does not fit into $X'$.  How much Riemannian geometry is
remembered by the symplectic geometry of the unit ball cotangent
bundle?

Our results give a pretty good understanding of this question
when $X$ and $X'$ are rectangles.  Let $X$ be the rectangle $[0,
L_1] \times ... \times [0, L_n]$ with the convention $L_1 \le ...
\le L_n$.  Let $X'$ be the rectangle $[0,L_1'] \times ... \times
[0, L_n']$ with $L_1' \le ... \le L_n'$.  The following
proposition from \cite{Gu1} describes when it is possible to find
an expanding embedding from $X$ into $X'$.

\begin{refprop} Up to a constant factor $C(n)$,
there is an expanding embedding from $X$ into $X'$ if and only if
$L_1 ... L_k \lesssim L_1' ... L_k'$ for each $k$ in the range $1
\le k \le n$.
\end{refprop}

Using Theorem 1, we get the following information about
symplectic embeddings of unit ball cotangent bundles.

\begin{refprop} Up to a constant factor $C(n)$ there is a
symplectic embedding from $U^* X$ into $U^* X'$ if and only $L_1
\lesssim L_1'$ and $L_1 ... L_n \lesssim L_1' ... L_n'$.
\end{refprop}

Hence the symplectic geometry remembers the two obstructions $L_1
\lesssim L_1'$ and $L_1 ... L_n \lesssim L_1' ... L_n'$ but
forgets the other n-2 obstructions.  In dimension $n=2$, the
symplectic geometry is roughly equivalent to the geometry of
expanding embeddings (at least for rectangles).  In dimension $n
\ge 3$, the symplectic geometry becomes more flexible than the
geometry of expanding embeddings.

\vskip5pt

{\bf A physical analogy}

\vskip5pt

Traynor's work proves Theorem 1 in the case $n=2$, so the
first new embedding happens in the case $n=3$.  Here is a typical
example.  Suppose that $P$ is equal to $B^2(\delta) \times B^2(1)
\times B^2(1)$ where $\delta > 0$ is an arbitrarily small number. 
Suppose that $P'$ is equal to $B^2(2 \delta) \times B^2(10
\delta) \times \mathbb{R}^2$.  In this paper, we will construct a
symplectic embedding from $P$ into $P'$.  Notice that in this
case, $R_1 R_2$ is much larger than $R_1' R_2'$.  I would like to
call this embedding the catalyst map because of the following
analogy with physics.

We think of a physical system with three degrees of freedom,
described by $(x_1,y_1)$, by $(x_2, y_2)$, and by $(x_3, y_3)$. 
For example, the system could consist of three particles each
moving in 1-dimensional space: $x_i$ would denote the position of
particle $i$, and $y_i$ would denote the momentum of particle
$i$.  At the initial time, we might know that the system lies in
the polydisk $P$.  In other words, we have detailed knowledge of
the state of particle 1, and medium knowledge of the states of
particles 2 and 3.  We would like to get more control over
particle 2, and we are willing to lose control of particle 3.  At
the end, we will probably discard particle 3 as exhaust.  The
system will evolve by a Hamiltonian diffeomorphism, and for the
purposes of this discussion let us suppose that we can apply to
the system any symplectic embedding.  If we just work with
particles 2 and 3, Gromov's non-squeezing theorem limits our
ability to do what we want to do: we cannot even symplectically
embed $B^2(1) \times B^2(1)$ into $B^2(1/2) \times \mathbb{R}^2$. 
We introduce a new, highly organized component, particle 1, and
then evolve the system by the catalyst map so that it lands in
$P'$.  By letting all three particles interact, we are able to
improve our knowledge of particle 2.  Particle 1 plays the role
of a catalyst: it comes out almost unchanged, but with its help,
particles 2 and 3 have interacted in a way they could not have
done on their own.  At the end of the interaction, our knowledge
of the catalyst has degraded by a factor of $2$, but we have
improved our knowledge of particle 2 by an arbitrary factor.

(It would be interesting to know if any real-world systems behave
in a way similar to the catalyst map.  There is an important
caveat about trying to relate the non-squeezing theorem to
practical problems in physics: a symplectic embedding can map an
arbitrary fraction of the volume of the unit ball into a thin
cylinder.  Given this caveat, I don't see any reason to think
physical systems would behave like the catalyst map.)

\vskip5pt

{\bf Outline}

\vskip5pt

We let $\Sigma$ denote a surface of genus 1 with one boundary
component, equipped with a symplectic form of area 1.  The key
step in our proof is the following lemma.

\begin{ML} For any radius $R$, the ball $B^4(R)$ symplectically
embeds into $\Sigma \times \mathbb{R}^2$.
\end{ML}

We prove the main lemma in Section 2.  The lemma builds on work
of Polterovich, who constructed a similar embedding using a
closed torus instead of $\Sigma$.  To prove the main lemma, we
modify Polterovich's construction in order to avoid a small disk
in the torus.

In Section 3, we use the main lemma to prove Theorem 1.  We also
give two corollaries.

\begin{introcor} For any radius $R$ and any $\epsilon > 0$, there
is an immersion from $B^4(R)$ into the standard cylinder $B^2(1)
\times \mathbb{R}^2$ so that each point in the range has at most
two preimages and so that the set of points with two preimages
has volume less than $\epsilon$.
\end{introcor}

This corollary shows that the non-squeezing theorem cannot be
weakened to allow immersions.  The other corollary of the main
lemma is a non-embedding result which builds on the non-squeezing
theorem.  Let $\Sigma(\epsilon)$ denote a surface with genus 1
and one boundary component equipped with a symplectic form of
area $\epsilon^2$.

\begin{introcor} Suppose that $\Sigma(\epsilon) \times B^2(R)$
symplectically embeds in the standard cylinder $B^2(1) \times
\mathbb{R}^2$.  No matter what the value of $\epsilon$, $R \le
1$.
\end{introcor}

In an appendix, we review the connection between expanding
embeddings and symplectic embeddings.  From this point of view,
we construct some symplectic embeddings of polydisks similar to
those of Traynor in \cite{T}.

\vskip5pt

{\bf Acknowledgements.} I would like to thank Helmut Hofer and Yasha
Eliashberg for helpful conversations.

\section{The main lemma}

Let $\Sigma$ denote a surface of genus 1 with one boundary
component.  We equip $\Sigma$ with a symplectic form of area 1. 
The key step in our proof is a symplectic embedding in four
dimensions described in the following lemma. We let $B^{2n}(R)$
denote the ball of radius $R$ in $\mathbb{R}^{2n}$ equipped with
the standard symplectic form.

\begin{ML} For any radius $R \ge 1/3$, there is a symplectic
embedding from $B^4(R)$ into $\Sigma \times B^2(10 R^2)$.
\end{ML}

This main lemma builds on a result of Polterovich.  Let $T^2(1)$
denote a torus equipped with a symplectic form of area 1.

\begin{reflemma} (Polterovich) For any radius $R \ge 1/3$, there is a
symplectic embedding from $B^4(R)$ into $T^2(1) \times B^2(10
R^2)$.
\end{reflemma}

Our lemma is a stronger version of Polterovich's.  We can
think of $\Sigma$ as $T^2(1)$ with a point $p$ removed.  To
prove our lemma, we have to modify Polterovich's embedding so
that its image avoids $\{ p \} \times \mathbb{R}^2$.

Before turning to the proof, let us indicate how to use the main
lemma.  Because of our lemma, we can embed $B^2(1/10) \times
B^4(R)$ into $B^2(1/10) \times \Sigma \times \mathbb{R}^2$. In
the next section, we will give a straightforward symplectic
embedding from the product $B^2(1/10) \times \Sigma$ into
$B^2(1/5) \times B^2(1)$.  Combining these two embeddings, we get
an embedding from $B^2(1/10) \times B^4(R)$ into $B^2(1/5) \times
B^2(1) \times \mathbb{R}^2$.  Up to scaling, this last embedding is
the catalyst map described in the introduction.  Notice that we
cannot use Polterovich's lemma in this argument because $T^2(1)$
does not symplectically embed into $\mathbb{R}^{4}$.

\vskip5pt

{\it Proof of Polterovich's lemma}: The embedding is given by a
linear map.  We write $B$ as an abbreviation of $B^4(R)$.  
Let $V$ denote an oriented 2-dimensional subspace of
$\mathbb{R}^4$ chosen so that the integral of $\omega$ over the
oriented disk $B \cap V$ is equal to $\pi/9$.  (We can find such a
subspace by continuity.)  Now there is a linear
symplectomorphism $L$ which maps planes parallel to $V$ to planes
parallel to the $x_1-y_1$ plane, and that maps disks parallel to
$V$ to disks.  The image $L(B)$ is an ellipsoid, and the
intersection of $L(B)$ with any plane parallel to the $x_1-y_1$
plane is equal to a disk of radius at most $1/3$.  The situation
is illustrated in Figure 1.

\includegraphics{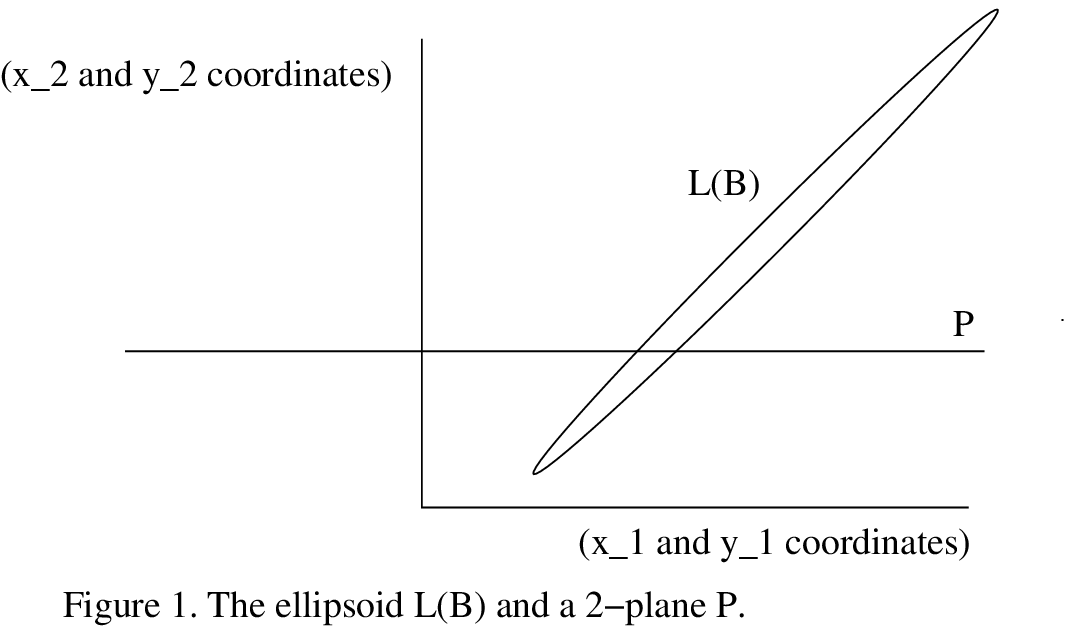}

(The plane $P$ in the figure is a 2-plane parallel to the $x_1-y_1$
plane.  In the figure, the plane $P$ meets $L(B)$ in a small
disk.  Nevertheless, the projection from $L(B)$ to the $x_1-y_1$
plane is large.)

We define $Q: \mathbb{R}^4 \rightarrow T^2(1) \times
\mathbb{R}^2$ to be the quotient map modding out by the integer
lattice generated by $x_1$ and $y_1$.  The quotient map $Q$ is
symplectic, and $Q$ restricted to $L(B)$ is an embedding.  To see
this, we just intersect $L(B)$ with any plane parallel to the
$x_1-y_1$ axis, and look at the quotient map from the
intersection to $T^2$.  Since the intersection is an open disk of
radius at most $1/3$, the quotient map is injective.  Hence $Q
\circ L$ is a symplectic embedding from $B(R)$ into $T^2(1)
\times \mathbb{R}^2$.

If we project $L(B)$ onto the $x_2-y_2$ plane, we get an ellipse
with area on the order of $R^4$.  We can choose the linear map
$L$ so that this ellipse is a disk.  Next we estimate the area of
this disk more carefully.  Let us denote its radius by $S$.  Let
$\pi_2$ denote the projection from $\mathbb{R}^4$ onto the
$x_2-y_2$ plane.  The point $\pi_2[L(0)]$ is the center $c$ of
the disk.  Hence the preimage $\pi_2^{-1}(c)$ is a disk of area
$\pi/9$.  Now if $p$ is any point in the concentric disk of half
the radius, then by convexity $\pi_2^{-1}(p)$ contains a disk of
area $\pi/36$.  Hence the volume of $L(B)$ is at least
$(\pi/36) \pi (S/2)^2 = (1/144) \pi^2 S^2$.  On the other
hand the volume of $L(B)$ is the same as the volume of $B$, which
is equal to $(1/2) \pi^2 R^4$.  Hence $S < (72)^{1/2} R^2
< 10 R^2$.  This finishes the proof of Polterovich's lemma.

\vskip5pt

\proof Now we turn to the proof of the main lemma.  We will
proceed by modifying Polterovich's embedding.  As before, we
write $B$ to abbreviate $B^4(R)$, and we let $L$ be the linear
map constructed above.

Now we outline our strategy.  We let $\Psi$ denote a
symplectomorphism of the $x_1-y_1$-plane which we will have to
choose carefully later on.  We write $\bar \Psi$ to denote the
product of $\Psi$ with the identity, which is a
symplectomorphism of $\mathbb{R}^4$.  We will follow the plan of
Polterovich's proof above except that we will use the non-linear
symplectomorphism $\bar \Psi \circ L$ in place of the linear
symplectomorphism $L$.  

Our embedding will be the composition $Q \circ \bar \Psi \circ
L$, which is automatically a symplectic immersion from $B^4(R)$
into $T^2(1) \times \mathbb{R}^2$.  We let $0$ denote the point
of $T^2(1)$ corresponding to the integer vectors in the $x_1-y_1$
plane, and we identify $\Sigma$ with $T^2(1) - \{ 0 \}$.  To
prove our lemma, we have to choose $\Psi$ so that $Q \circ \bar
\Psi \circ L (B)$ lands inside of $\Sigma \times \mathbb{R}^2$,
and so that $Q$ restricted to $\bar \Psi \circ L(B)$ is an
embedding.

We will choose the map $\Psi$ to obey the two properties below.

We let $\pi$ denote the projection from $\mathbb{R}^4$ onto the
$x_1-y_1$ plane.  We choose a number $\rho$ large enough so that
$\pi(L(B))$ lies in the disk of radius $\rho$ around the origin.

Property 1. The map $\Psi$ takes the disk of radius $\rho$
into the complement of all integer lattice points.

Let $S$ be a subset of the $x_1-y_1$ plane.  We say that $S$ is
aperiodic if the difference of any two points in $S$ is never a
non-trivial integer vector.

Property 2. If $D$ denotes any disk in the $x_1-y_1$ plane of
radius $1/3$, then $\Psi(D)$ is aperiodic.

Now we check that if $\Psi$ obeys Properties 1 and 2, then $Q
\circ \bar \Psi \circ L$ embeds $B^4(R)$ into $\Sigma \times
B^2(10 R^2)$.  First of all, we have to check that $Q \circ \bar
\Psi \circ L(B)$ lands inside of $\Sigma \times B^2(10 R^2)$.  We
know that $L(B)$ lands in $B^2(\rho) \times B^2(10 R^2)$.  By
Property 1, we see that $\bar \Psi \circ L(B)$ lands in
$(\mathbb{R}^2 - \mathbb{Z}^2) \times B^2(10 R^2)$.  And so $Q
\circ \bar \Psi \circ L(B)$ lands in $\Sigma \times B^2(10 R^2)$. 

Second, we have to check that $Q \circ \bar \Psi \circ L$ is an
embedding from $B$.  Let $p,q$ be points in $B$, and
suppose that $Q \circ \bar \Psi \circ L(p) = Q \circ \bar \Psi
\circ L(q)$.  It follows that $\bar \Psi \circ L(p) = \bar \Psi
\circ L(q) + (m,n,0,0)$ for some integers $m$ and $n$.  In
particular we see that $\bar \Psi \circ L(p)$ and $\bar \Psi
\circ L(q)$ have the same $x_2,y_2$ coordinates. Since $\bar
\Psi$ doesn't change $x_2,y_2$ coordinates, it follows that
$L(p)$ and $L(q)$ have the same $x_2, y_2$ coordinates.  Let $W$
denote the 2-plane of all points in $\mathbb{R}^4$ with the same
$x_2$ and $y_2$ coordinates as $L(p)$ and $L(q)$.  The
intersection $W \cap L(B)$ lies in a disk $D$ of radius $1/3$. 
In particular $L(p)$ and $L(q)$ both lie in $D$.  We can think of
$x_1$ and $y_1$ as coordinates on this plane $W$, and so we can
define $\Psi(D) \subset W$.  Now $\bar \Psi \circ L(p)$ and $\bar
\Psi \circ L(q)$ both lie in $\Psi(D)$.  By Property 2, we know
that $\Psi(D)$ is aperiodic.  On the other hand, we established
above that $\bar \Psi \circ L(p) = \bar \Psi \circ L(q) + (m,n,
0,0)$.  Therefore, $m$ and $n$ are zero, $p$ is equal to $q$, and
$Q \circ \bar \Psi \circ L$ does embed $B$ into $\Sigma \times
B^2(10 R^2)$ as claimed.

It remains to construct the map $\Psi$ with the two properties
above.  As a tool for constructing $\Psi$, we define a
diffeomorphism $\Phi$ of the $x_1$-axis.  The diffeomorphism
$\Phi$ maps each integer point to itself.  It is periodic with
period 1.  The derivative $d \Phi$ is always at least $9/10$.  At
each integer point $d \Phi = 100 \rho$.  Finally, the
displacement $|\Phi(x) - x|$ is at most $10^{-4}$ for every $x$. 
Such a diffeomorphism is easy to find.

Now we define $\Psi$ by the following formula.

$$\Psi(x, y) = (\Phi(x), (1/2) + [d \Phi(x)]^{-1} y).$$

To give some sense of this map, we sketch the image of the ball
of radius 3 under the map $\Psi$.

\vskip15pt

\includegraphics{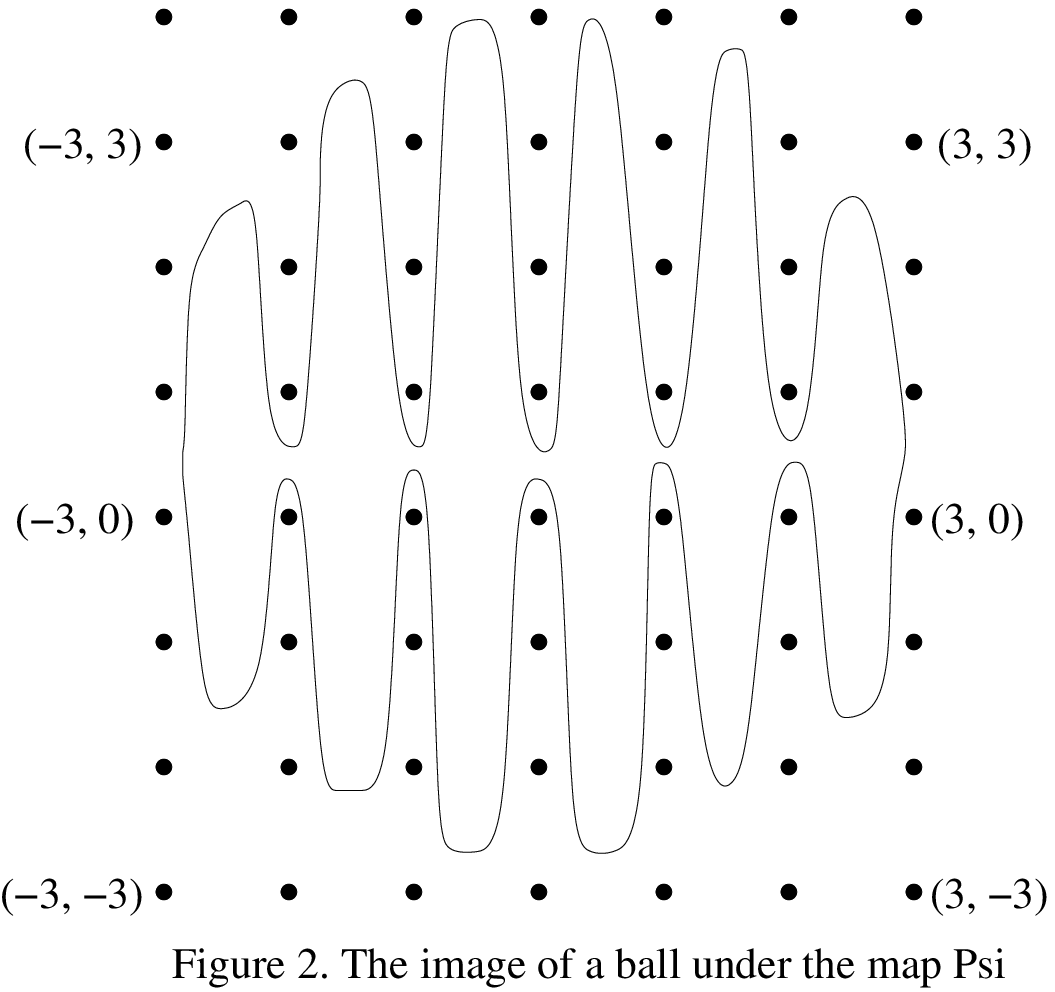}

\vskip10pt

The map $\Psi$ is clearly a diffeomorphism.  Next we check that
it is area preserving by computing its Jacobian.  We let
$\Psi_1$ be the x-coordinate of $\Psi$ and $\Psi_2$ be the
y-coordinate.  We compute some derivatives

$\partial \Psi_1 / \partial x = d \Phi(x).$

$\partial \Psi_1 / \partial y = 0.$

$\partial \Psi_2 / \partial y = [d \Phi(x)]^{-1}.$

It follows that the determinant of $d \Psi$ is equal to 1.  Hence
$\Psi$ is a symplectomorphism.

Now we check that $\Psi$ obeys Properties 1 and 2.

Proof of Property 1.  Suppose that $(x,y)$ is in the disk of
radius $\rho$.  We have to check that $\Psi(x,y)$ is not an
integer lattice point.  Suppose that $\Psi(x,y)$ has
$x$-coordinate equal to an integer.  In this case, $\Phi(x)$ is
an integer.  By the definition of $\Phi$ it follows that $x$ is
an integer, and so $d \Phi(x) = 100 \rho$.  But the norm $|y|$ is
at most $\rho$.  Therefore, the y-coordinate of $\Psi(x,y)$ is
between $(1/2) - (1/100)$ and $(1/2) + (1/100)$.  In particular,
the y-coordinate is not an integer.

Proof of Property 2. Suppose that $D$ is a disk of radius $1/3$. 
We have to check that $\Psi(D)$ is aperiodic.  Let $p$ and $q$
be two points in $D$.  We have to check that $\Psi(p) -
\Psi(q)$ is not a non-trivial integer lattice point.  We suppose
that $\Psi(p) - \Psi(q)$ is an integer lattice point and prove
that $p=q$.

From the definition of $\Phi$, we know that the displacement
$|\Phi(x) - x|$ is at most $10^{-4}$.  Therefore, the
x-coordinate of $\Psi(p)$ agrees with the x-coordinate of $p$ up
to an error of $10^{-4}$.  Similarly for q.  Therefore,
$|\Psi_1(p) - \Psi_1(q)| \le (2/3) + 2 \cdot 10^{-4}$.  Since the
difference $\Psi_1(p) - \Psi_1(q)$ is an integer, we conclude
that the difference is zero.  But the x-coordinate $\Psi_1(p)$
only depends on the x-coordinate of $p$, and so the x-coordinates
of $p$ and $q$ are equal.  Let $x_0$ be the x-coordinate of each
point.

Now the y-coordinate of $\Psi(p)$ is $1/2 + [d \Phi(x_0)]^{-1}
y(p)$ and the y-coordinate of $\Psi(q)$ is $1/2 + [d
\Phi(x_0)]^{-1} y(q)$.  Hence their difference is $[d
\Phi(x_0)]^{-1} (y(p) - y(q))$.  But by the definition of $\Phi$,
$d \Phi$ is at least $9/10$.  Therefore, the difference has norm
at most $(10/9) |y(p) - y(q)| \le (10/9) (2/3) < 1$.  Since we
assumed the difference is an integer, the difference is zero. 
Then it follows that $y(p) = y(q)$ and finally that $p=q$.  This
finishes the proof of Property 2 and hence the proof of the main
lemma.  \endproof

\section{Consequences of the main lemma}

In this section, we use the main lemma to prove Theorem 1 and
afterwards give some other minor consequences.  To prove Theorem
1 we need one more lemma, which gives us a way to embed $B^2(W)
\times \Sigma$ into a polydisk.

\begin{lemma} If $W \le 1/10$, then $B^2(W) \times \Sigma$
symplectically embeds in $B^2(2W) \times B^2(1)$.
\end{lemma}

\proof We begin by choosing a symplectic immersion $I$ of
$\Sigma$ into $B^2(1)$.  We choose the immersion so that the
image of $\Sigma$ meets the unit square $[-1/2, 1/2]^2$ in a
particularly simple form.  Namely, the image contains the strip
$S = [-1/2, 1/2] \times [-W, W]$ and the strip $S' = [- W,
W]
\times [-1/2, 1/2]$.  Other than these two strips, the immersion
does not hit the square $[-1/2, 1/2]^2$.  We can arrange that the
immersion $I$ is an embedding except for the overlap of these two
strips.  The immersion is illustrated in Figure 3.

\includegraphics{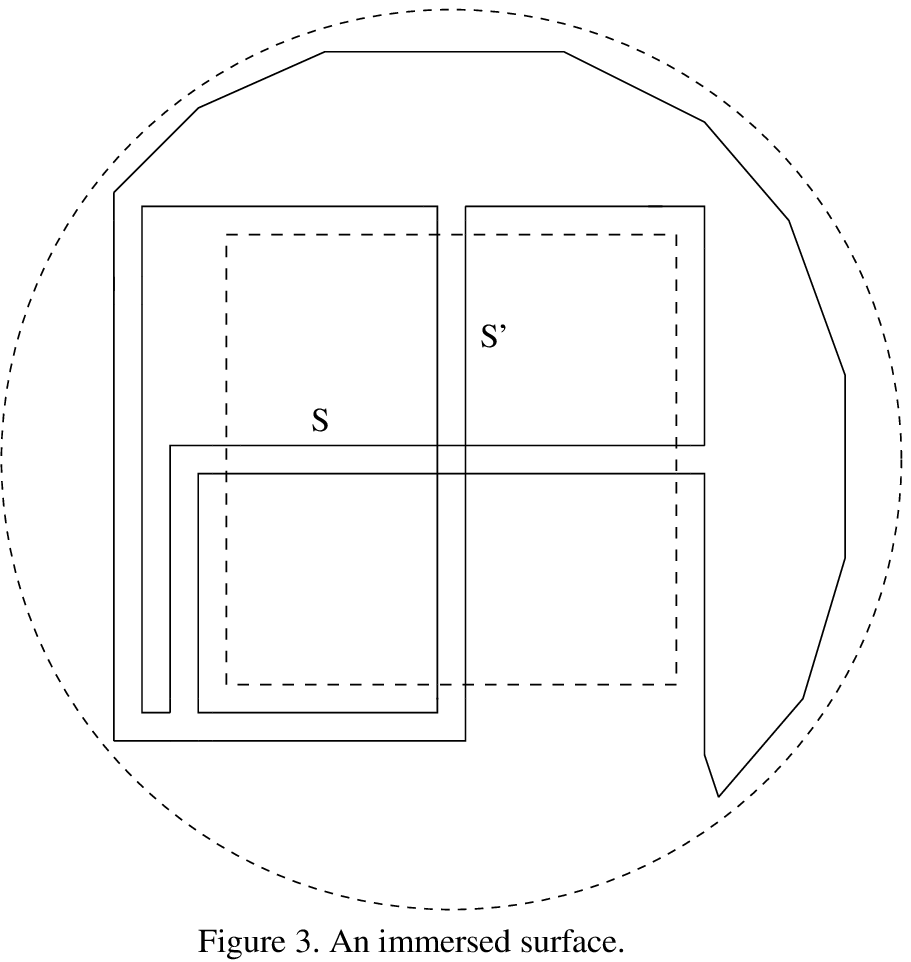}

The dotted circle on the outside of the figure is the unit
circle.  The dotted square is the unit square $[-1/2, 1/2]^2$. 
The solid lines are the boundary of the immersed copy of
$\Sigma$.  The horizontal strip $S$ and the vertical strip $S'$
are labelled in the figure.

Now the map $I' = I \times id$ gives a symplectic immersion from
$\Sigma \times B^2(W)$ into $B^2(1) \times \mathbb{R}^2$.  We
will modify this immersion to make it a sympectic embedding. 
Roughly, we are going to lift the strip $S \times B^2(W)$ in the
$y_2$ direction a distance $2W$ in order to push it just over the
strip $S' \times B^2(W)$.

We will modify $I'$ on the region $S \times B^2(W)$ using a
Hamiltonian flow.  We use coordinates $x_1, y_1$ on $B^2(1)$ and
coordinates $x_2, y_2$ on $\mathbb{R}^2$.  We use the Hamiltonian
$H = \Psi(x_1) x_2$, where $\Psi$ is a bump function, equal to
$1$ on $[-1/6, 1/6]$, equal to $0$ outside of $[-1/3, 1/3]$, and
with $|\nabla \Psi| \le 7$.  We run the Hamiltonian flow for time
$T=2W$.  This defines a family of symplectic embeddings of $S
\times B^2(W)$ into $B^2(1) \times
\mathbb{R}^2$.  All of these embeddings agree with $I'$ on a
neighborhood of $\pm 1/2 \times [- W, W] \times
B^2(W)$, so they extend to immersions of $\Sigma \times
B^2(W)$ into $B^2(1) \times \mathbb{R}^2$.

We claim that at time $T = 2W$ we get an embedding of $\Sigma
\times B^2(W)$ into $B^2(1) \times \mathbb{R}^2$.  We used the
Hamiltonian flow above to change the embedding of $S
\times B^2(W)$.  Meanwhile, the
embedding of $(\Sigma - S) \times B^2(W)$ is left unchanged.  We
now have to check that our new embedding of $S \times B^2(W)$ is
disjoint from the original embedding of $(\Sigma - S) \times
B^2(W)$.

During this calculation it's helpful to notice that the
Hamiltonian flow leaves $x_1$ and $x_2$ invariant.  As a first
step, we check that at time $T=2W$, our embedding of $S \times
B^2(W)$ lands in $[-1/2, 1/2]^2 \times \mathbb{R}^2$.  Since
$x_1$ doesn't change during the flow, it follows that the
$x_1$-coordinate of our embedding of $S \times B^2(W)$ lies in
$[-1/2, 1/2]$.  Next we deal with the $y_1$-coordinate.  At the
initial time $|y_1| \le W$.  During the Hamiltonian flow, $y_1$
changes at the rate $\Psi'(x_1) x_2$.  The gradient $\Psi'(x_1)$
has norm at most 7.  Also, the $x_2$-coordinate doesn't change
during the flow and it has norm at most $W$.  Therefore, $y_1$
changes at a rate at most $7 W$.  Hence, at time $T = 2W$, $|y_1|
\le W + 7 W (2 W)$.  Because $W \le 1/10$, we conclude that at
the final time $|y_1| < 1/2$. 

The only part of $(\Sigma - S) \times B^2(W)$ that lies in
$[-1/2, 1/2]^2 \times \mathbb{R}^2$ is the other strip $S'
\times B^2(W)$.  We check that the image of $S \times B^2(W)$ is
disjoint from the vertical strip $S' \times B^2(W)$.  Suppose
that $(x_1, y_1, x_2, y_2)$ lies in the image of our embedding of
$S \times B^2(W)$.  If $|x_1| > 1/6 > W$, then this point does
not lie in $S' \times B^2(W)$.  It remains to consider the case
that $|x_1| \le 1/6$.  We recall that $x_1$ did not change during
the Hamiltonian flow.  The derivative of $y_2$ during the flow
was $\Psi(x_1) = 1$.  Since we ran the flow for time $2W$, the
value of $y_2$ increased during the flow by $2W$.  Since $y_2$
was initially in $(-W, W)$, at time $T=2W$ we have $y_2 > W$. 
Hence our point is disjoint from $S' \times B^2(W)$.  This
argument shows that we have defined an embedding from $\Sigma
\times B^2(2W)$ into $B^2(1) \times \mathbb{R}^2$.

The final two coordinates of our embedding always have the form
of a point in $B^2(W)$ plus a vector in the positive $y_2$
direction of length at most $2W$.  Therefore the image lies in
$B^2(1) \times B^2[(0,W), 2W]$, where the second term denotes the
ball of radius $2W$ around the point $(0,W)$ in the $x_2-y_2$
plane.
\endproof

Combining this lemma with the main lemma we can construct the
catalyst map described in the introduction.  Let $R$ denote a
large radius.  By the main lemma, we can embed $B^2(1/10) \times
B^2(R) \times B^2(R)$ into $B^2(1/10) \times \Sigma \times B^2(20
R^2)$.  Now applying the last lemma, we can embed this shape into
$B^2(1/5) \times B^2(1) \times B^2(20 R^2)$.  After scaling, this
map is the catalyst map.  Note that the first radius gets bigger
by a factor of only 2, while the second radius gets smaller by a
factor of $R$.

We remark that this construction is not completely explicit.  The
surface $\Sigma$ appears in two different ways.  First, we take
the $x_1 - y_1$ plane, mod out by the action of $\mathbb{Z}^2$,
and then remove a small disk or a point.  Second, we take the
immersed surface in Figure 3 above with the induced symplectic
form.  These two surfaces are symplectomorphic by Moser's theorem
\cite{M}.  The catalyst map is the composition of three steps:
the map from the main lemma, then a Moser symplectomorphism, and
then the map constructed in Lemma 3.1.

To get further embeddings, we need to combine our construction
with embeddings coming from Traynor's work.  Also, to help keep
track of what embeds in what, we adopt the following notation. 
If $P$ is a polydisk with radii $R_1 \le ... \le R_n$, and if $C
> 0$, then we write $C P$ to denote the magnified
polydisk with radii $C R_1 \le ... \le C R_n$.

\begin{prop} (Traynor) There exists a constant $C > 0$ so that the
following holds.  Suppose $P$ is a polydisk with radii $R_1 \le
R_2$.  Suppose $1 \le \lambda \le (R_2/ R_1)^{1/2}$.  Let $P'$ be
the polydisk with $R_1' = \lambda R_1$, and $R_2' = R_2 /
\lambda$.  Then $P$ symplectically embeds into $C P'$.
\end{prop}

(If $R_1' = R_2'$, then this result follows immediately from
Theorem 1.3 of \cite{T}.  The general case can be proved using
the same method.  In the appendix, we sketch a proof of this
proposition using expanding embeddings.)

By combining these tools, we get the following generalization of
the catalyst map.

\begin{prop} There exists a constant $C > 0$ so that the
following holds.  Suppose $P$ is a polydisk with radii $R_1 \le
R_2 \le R_3$.  Suppose $1 \le \lambda \le R_2/ R_1$.  Let $P'$ be the
polydisk with $R_1' = R_1$, $R_2' = R_2 / \lambda$ and
$R_3' = R_3 \lambda$.  Then $P$ symplectically embeds into
$C P'$.
\end{prop}

\proof By Proposition 1 (due to Traynor), we know that $B^2(R_2)
\times B^2(R_3)$ symplectically embeds in $B^4(C R_2^{1/2}
R_3^{1/2})$.  Applying the main lemma and a scaling argument,
this ball symplectically embeds in $\Sigma(R_2') \times B^2(C
R_3')$, where $\Sigma(R_2')$ denotes a surface of genus 1 with one
boundary component equipped with a symplectic form of area
$(R_2')^2$.  Hence $B^2(R_1) \times B^2(R_2) \times B^2(R_3)$
symplectically embeds into $B^2(R_1') \times \Sigma(R_2') \times
B^2(C R_3')$.  By the last lemma, $B^2(R_1') \times \Sigma(R_2')$
symplectically embeds in $B^2(C R_1') \times B^2(C R_2')$.  Hence
$P$ symplectically embeds in $C P'$. \endproof

Our main theorem follows from combining these two propositions. 
It requires no new ideas, but the algebra is a bit tedious.

\begin{theorem} For each integer $n$, there is a constant $C(n)$
so that the following holds.  Let $P$ and $P'$ be polydisks of
dimension $2n$.  Suppose that $R_1 \le R_1'$ and $R_1 ... R_n \le
R_1' ... R_n'$.  Then $P$ symplectically embeds in $C(n) P'$.
\end{theorem}

\proof Using Proposition 1 repeatedly, we see that $P$ embeds
symplectically in $C(n) P(1)$, where $R(1)_1 = R_1$ and $R(1)_i = (R_2
... R_n)^{\frac{1}{n-1}}$ for $2 \le i \le n$.  Again, using
Proposition 1 repeatedly, we see that $P(1)$ embeds
symplectically in $C(n) P(2)$, where $R(2)_1 = R_1'$ and $R(2)_i
= [R_1 ... R_n / R_1']^{\frac{1}{n-1}}$.  Finally, using
Proposition 2 repeatedly, we see that $P(2)$ embeds
symplectically in $C(n) P'$. \endproof

To end this section, we give two more consequences of the
main lemma.

\begin{cor} For any $R$ and any $\epsilon$, there is a
symplectic immersion from $B^4(R)$ into $B^2(1) \times
\mathbb{R}^2$ so that each point in the range has
at most two preimages and so the set of points with two preimages
has volume at most $\epsilon$.
\end{cor}

\proof First we use the main lemma to embed $B^4(R)$ into $\Sigma
\times B^2(10 R^2)$.  Next, we pick a symplectic immersion from $\Sigma$
into $B^2(1)$.  We can choose this immersion so that each point
in the target has at most two preimages and so that the set of
double points has area at most $\delta$ for any $\delta > 0$. 
(The immersion we need is illustrated in Figure 3 above.) 
Composing the immersion and the embedding we get a symplectic
immersion from $B^4(R)$ into $B^2(1) \times B^2(10 R^2)$.  It has
at worst double points and the set of double points has area at
most $10 \delta R^2$. \endproof

Lastly we give a non-embedding result.  Let $\Sigma(\epsilon)$
denote a rescaling of $\Sigma$ with symplectic area $\epsilon^2$.

\begin{cor} If $\Sigma(\epsilon) \times B^2(W)$ symplectically embeds in
the cylinder $B^2(1) \times \mathbb{R}^2$ then $W \le 1$,
regardless of $\epsilon$. 
\end{cor}

\proof  Suppose $\Sigma \times B^2(W)$ symplectically embeds
in $B^2(R) \times \mathbb{R}^2$.  We know that $B^4(W)$ embeds in
$\Sigma \times \mathbb{R}^2$, and so $B^6(W)$ embeds in $\Sigma
\times B^2(W) \times \mathbb{R}^2$ which embeds in $B^2(R)
\times \mathbb{R}^4$.  By Gromov's non-squeezing
theorem, we conclude that $R \ge W$.

Now if $\Sigma(\epsilon) \times
B^2(W)$ symplectically embeds in $B^2(1) \times \mathbb{R}^2$,
we can scale the domain and range to symplectically embed $\Sigma
\times B^2(\epsilon^{-1} W)$ into $B^2(\epsilon^{-1}) \times
\mathbb{R}^2$.  By the last paragraph, we get $\epsilon^{-1} W
\le \epsilon^{-1}$, and so $W \le 1$. \endproof

I don't know whether this result is sharp.  Using Lemma 3.1., we
can symplectically embed $\Sigma(\epsilon) \times B^2(W)$ into
$B^2(2 W) \times
\mathbb{R}^2$ for any values of $\epsilon$ and $W$.  The
construction in Lemma 3.1 can be improved by using a square in
the $x_2-y_2$ coordinates instead of a disk.  With this
improvement, we can symplectically embed $\Sigma(\epsilon) \times
B^2(W)$ into $B^2(\sqrt2 W) \times \mathbb{R}^2$.  I don't know
whether the factor $\sqrt 2$ can ever be reduced.

\section{Appendix: Expanding embeddings and symplectic
embeddings}

In this section we sketch a proof of Proposition 1 using
expanding embeddings of rectangles. 

If $M$ is a smooth manifold, then the cotangent bundle $T^* M$
has a canonical symplectic structure.  If $M$ is a Riemannian
manifold, then we let $U^* M$ denote the unit ball cotangent
bundle.  If $M, N$ are Riemannian manifolds, then an embedding
$I$ from $M$ to $N$ is called expanding if for any tangent vector
$v$ in $TM$, $|dI(v)| \ge |v|$.  An expanding embedding increases
the length of every curve.  Any embedding from $M$ to $N$ induces
a symplectic embedding from $T^* M$ into $T^* N$.  An expanding
embedding from $M$ to $N$ induces a symplectic embedding from
$U^*M$ to $U^*N$.

Suppose that $X$ and $X'$ are two-dimensional rectangles: $X =
[0, L_1] \times [0, L_2]$, with $L_1 \le L_2$ and $X' = [0, L_1']
\times [0, L_2']$ with $L_1' \le L_2'$.  If $L_1 \le L_1'$
and $L_1 L_2 \le L_1' L_2'$, then there is an expanding embedding
from $X$ into $5 X'$.  This embedding is illustrated in Figure 4.

\includegraphics{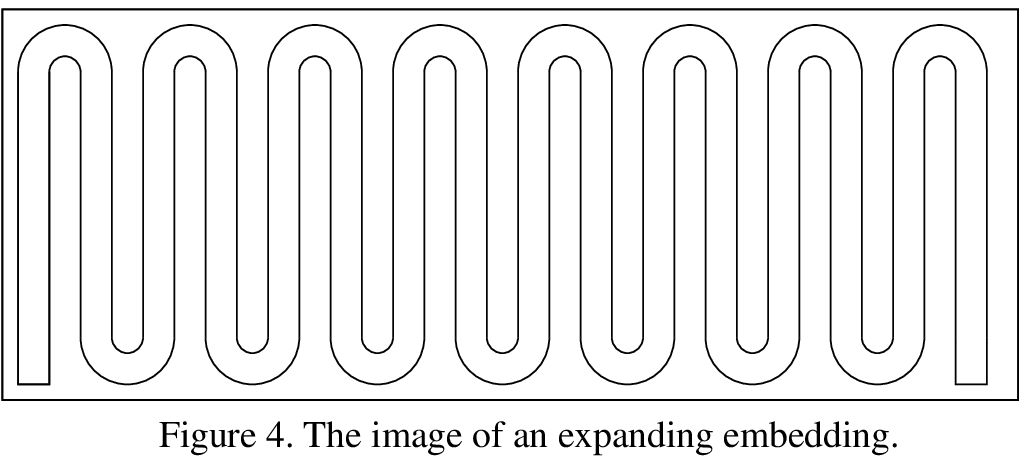}

The thicker rectangle on the outside is $5 X'$.  The snake-like
shape inside is the image of $X$ by an expanding embedding.

This expanding embedding induces a symplectic embedding from
$U^* X$ to $U^* (5 X')$.

The unit ball cotangent bundle $U^* X$ is not a polydisk, but it's
close enough to a polydisk to prove Proposition 1.  In
particular, $U^* X$ contains the set of
all pairs $(x,v)$, where $x$ is a point in $X$ and $v$ is a
cotangent vector $(v_1, v_2)$ with $|v_1| \le 2^{-1/2}$ and
$|v_2| \le 2^{-1/2}$.  The latter shape is a 4-dimensional
rectangle and so it is symplectomorphic to a polydisk $P =
B^2(R_1) \times B^2(R_2)$ with radii satisfying $\pi
R_1^2 = 2^{-1/2} L_1$ and $\pi R_2^2 = 2^{-1/2} L_2$.

Similarly, $U^* (5 X')$ is contained in the
rectangle of points $(x,v)$, where $x$ is a point in $5 X'$ and
$v = (v_1, v_2)$ is a vector with $|v_i| \le 1$.  This rectangle
is symplectomorphic to a polydisk $P' = B^2(R'_1) \times
B^2(R'_2)$, with radii satisfying $\pi (R'_1)^2 = 5 L_1'$ and $\pi
(R'_2)^2 = 5 L_2'$.

To summarize, $P$ symplectically embeds into $U^* X$.  Then
$U^*X$ symplectically embeds into $U^* (5 X')$.  Finally, $U^* (5
X')$ symplectically embeds into $P'$.  Hence $P$ symplectically
embeds in $P'$.  If we begin with $P$ and $P'$ obeying $3 R_1 \le
R_1'$ and $9 R_1 R_2 \le R_1' R_2'$, then a calculation shows
that we can symplectically embed $P$ into $P'$ using this
construction by choosing $X$ and $X'$ appropriately.  Hence, this
construction gives a proof of Proposition 1 with constant $C =
3$.

\end{document}